\documentclass[a4paper,leqno,10pt]{article}
\usepackage{amssymb}
\usepackage{amsmath}
\usepackage[german]{babel}
\begin{document}
\def\refname{ }
\begin{center}
\begin{Large}
\textbf{On Pad\'e approximations and global relations of some Euler-type series}
\end{Large}
\vskip0.4cm
\sc{Keijo V\"a\"an\"anen}
\end{center}
\vskip0.6cm

\begin{center}
Abstract
\end{center}

We shall consider some special generalizations of Euler's factorial series. First we construct Pad\'e approximations of the second kind for these series. Then these approximations are applied to study global relations of certain $p$-adic values of the series.\\

\noindent
2010 Mathematics Subject Classification: 11J13 (Primary), 11J61, 11J72 (Secondary)

\noindent
Keywords: $p$-adic linear form, Euler-type function, global relation, Pad\'e approximation

\section{Introduction}

In a recent work \cite{MZ} Matala-aho and Zudilin study global relations of the famous Euler's factorial series
\begin{equation}\label{1}
E(z) = \sum_{n=0}^\infty n!(-z)^n,
\end{equation}
converging $p$-adically in the disc $\left|z\right|_p \leq 1$ for all primes $p$. They introduce Pad\'e approximations of the hypergeometric series
\begin{equation}\label{2}
_2F_0(\alpha,1\mid z) = \sum_{n=0}^\infty(\alpha)_nz^n,
\end{equation}
where the Pochhammer notation $(\alpha)_n$ is defined by $(\alpha)_0 = 1, (\alpha)_n = \alpha(\alpha+1)\cdots(\alpha+n-1), n \geq 1$, and use these approximations with $\alpha = 1$ to prove the following result, where $E_p(z)$ denotes the series $E(z)$ considered as a function in $\mathbb{Q}_p$, the $p$-adic completion of $\mathbb{Q}$.\\

\noindent
\textbf{Theorem} (Matala-aho, Zudilin). \textit{Given $a \in \mathbb{Z}\setminus\{0\}$, let $P_1$ be a subset of prime numbers such that
\[
\limsup_{n\rightarrow\infty}c^n n!\prod_{p\in P_1}\left|n!\right|_p^2 = 0, \ c= 4\left|a\right|\prod_{p\in P_1}\left|a\right|_p^2.
\]
Then either there exists a prime $p\in P_1$ for which $E_p(a)$ is irrational, or there are two distinct primes $p, q \in P_1$ such that $E_p(a) \neq E_q(a)$ (while $E_p(a), E_q(a) \in \mathbb{Q}$).}\\

For more information on the results and open questions concerning Euler's factorial series we refer to \cite{C2} and \cite{MZ}.

In the present work we shall consider several series
\begin{equation}\label{3}
\varphi_j(z) = _2F_0(\alpha_j,1\mid z) = \sum_{n=0}^\infty(\alpha_j)_nz^n, \ j = 1,\ldots,m,
\end{equation}
where $\alpha_1,\ldots,\alpha_m \neq 0,-1,-2,\ldots$ are rational numbers such that $\alpha_i-\alpha_j \notin \mathbb{Z}$, if $i \neq j$. If $\alpha_j = r_j/s_j$, where $r_j$ and $s_j \geq 1$ are coprime integers, then the series $\varphi_j(z)$ converges in $\mathbb{Q}_p$ for all $\left|z\right|_p \leq 1$, if $p \nmid s_j$, and for all $\left|z\right|_p < p^{-ord_p(s_j)}$, if $p \mid s_j$. Therefore a linear form
\[
L = \ell_0 + \ell_1\varphi_1(a) +\cdots+ \ell_m\varphi_m(a)
\]
is defined in $\mathbb{Q}_p$ for all $a\in \mathbb{Z}\setminus\{0\}, (\ell_0,\ell_1,\ldots,\ell_m)\in \mathbb{Z}^{m+1}\setminus\{\underline{0}\}$ and $p\in P_1$, if $P_1$ is a subset of prime numbers such that it does not contain any of the prime factors of $s:=lcm\{s_1,\ldots,s_m\}$. We shall use the notation $L_p$ to denote $L$ as an element of $\mathbb{Q}_p$ and call a relation $L=0 \ P_1$-global, if $L_p=0$ for all $p\in P_1$. Our main purpose here is to construct explicitly Pad\'e approximations of the second kind for the series $\varphi_j(z)$ and use these to prove the following result on $P_1$-global relations.\\

\noindent
\textbf{Theorem 1}. \textit{Let $P_1$ denote a subset of primes such that $p\notin P_1$, if $p\mid s$, and let $a\in \mathbb{Z}\setminus\{0\}$. There exists a positive constant $c$ (given explicitly in (17)) depending on $\alpha_1,\ldots,\alpha_m$ and $m$ such that there are no $P_1$-global relations if}
\begin{equation}\label{4}
\limsup_{n\rightarrow\infty}c^n\left|a\right|^{mn+1}(mn)!\prod_{p\in P_1}\left|a^{mn+n+1}(mn)!(n+1)!\right|_p = 0.\\
\end{equation}

Clearly the above condition (\ref{4}) holds if the complement of $P_1$ in the set of all primes is finite. It is also possible to study the number of primes $p$ with $L_p \neq 0$ or to give bounds in terms of $h:=\max\{\left|\ell_0\right|,\left|\ell_1\right|,\ldots,\left|\ell_m\right|\}$ for a prime $p$ with $L_p\neq 0$ along the lines of the papers \cite{BCY} and \cite{C2} considering more general classes of Euler-type series. Our series $\varphi_j(z)$ belong to the class studied in \cite{BCY}. Concerning the effectivity of the bounds the main difficulty in \cite{BCY} comes from the constant in Shidlovskii's lemma needed for the proof of the non-vanishing of the crucial determinant corresponding somehow our $\Omega$ in Lemma 1 below. Because of explicit construction we avoid here such difficulties. As an example of results obtained by our approximations we give the following theorem.\\
 
\noindent
\textbf{Theorem 2}. \textit{There are infinitely many primes $p$ satisfying $L_p \neq 0$. Moreover, we have $L_p \neq 0$ for some $p$},
\[
e^{\sqrt{\log \tilde{H}}} < p \leq m\tilde{H},  \quad \tilde{H} = \left[\frac{\log H}{\log\log H}(1+\frac{\frac{5}{4}(m+1)+6}{\sqrt{\log\log H}})\right],
\]
\textit{where $H = \max \{h,H_0\}$ with $h=\max\{\left|\ell_0\right|,\left|\ell_1\right|,\ldots,\left|\ell_m\right|\}$ and a positive constant $H_0$ \ (given explicitly in (\ref{20})) depending on $\alpha_1,\ldots,\alpha_m, m$ and $a$.}\\ 

In section 2 below we shall give explicit constructions of Pad\'e approximations of the second kind for the series $\varphi_j(z)$. More precisely, for given positive integers $n_1,\ldots,n_m, n \geq \max \{n_1,\ldots,n_m\}$, this is a system of polynomials $Q(z) \neq 0, P_1(z),\ldots, P_m(z)$ such that
\begin{equation}\label{5}
\deg Q(z) \leq N:= n_1+\cdots+n_m, \deg P_j(z) \leq N+n-n_j, ord_{z=0}(Q(z)\varphi_j(z)-P_j(z)) \geq N+n+1.
\end{equation}
Here we shall obtain such polynomials by using a simple method from \cite{F} and \cite{V}. After that the denominators of the coefficients of these polynomials are studied in section 3 and upper bounds for $\left|Q(z)\right|, \left|P_1(z)\right|,\ldots, \left|P_m(z)\right|$ and also for $\left|R_j(z)\right|_p, p\in P_1$, where $R_j(z) = Q(z)\varphi_j(z)-P_j(z)$, are obtained in section 4. Theorems 1 and 2 are then proved in the final section.

\section{Pad\'e approximations}

In this section we consider $\varphi_j(z)$ as formal power series in $\mathbb{Q}[[z]]$. For the construction of polynomials satisfying (\ref{5}) we denote
\[
Q(z) = \sum_{k=0}^N a_kz^k.
\]
Then
\[
Q(z)\varphi_j(z) = \sum_{\mu=0}^\infty c_{j\mu}z^\mu, \quad c_{j\mu} = \sum_{k=0}^{\min\{\mu,N\}} a_k(\alpha_j)_{\mu-k}, \ j=1,\ldots,m.
\] 
If we now choose the coefficients $a_k$ in such a way that at least one $a_k \neq 0$ and $c_{j\mu} = 0 \ (\mu = N+n-n_j+1,\ldots,N+n; j=1,\ldots,m)$, then the polynomials $Q(z)$ and
\[
P_j(z) = \sum_{\mu=0}^{N+n-n_j} c_{j\mu}z^\mu, \ j=1,\ldots,m,
\]
give the desired approximations. Thus we need to find a non-trivial solution to the system
\[
a_0(\alpha_j)_\mu+a_1(\alpha_j)_{\mu-1}+\cdots+a_N(\alpha_j)_{\mu-N} = 0, \ \mu=N+n-n_j+1,\ldots,N+n; j=1,\ldots,m,
\]
of linear homogeneous equations.

By denoting $a_{N-k} = b_k \ (k=0,1,\ldots,N)$ and
\[
\gamma_1 = n-n_1+1+\alpha_1, \ldots, \gamma_{n_1} = n+\alpha_1,
\]
\[
\gamma_{n_1+1} = n-n_2+1+\alpha_2, \ldots, \gamma_{n_1+n_2} = n+\alpha_2, \ldots
\]
\[
\gamma_{n_1+\cdots+n_{m-1}+1} = n-n_m+1+\alpha_m, \ldots, \gamma_N = n+\alpha_m,
\]
the above system can be given in the form
\begin{equation}\label{6}
b_0+b_1\gamma_i+b_2\gamma_i(\gamma_i+1)+\cdots+b_{N-1}\gamma_i(\gamma_i+1)\cdots(\gamma_i+N-2)=
\end{equation}
\[
-b_N\gamma_i(\gamma_i+1)\cdots(\gamma_i+N-1) , \ i = 1,\ldots,N.
\]
The coefficient determinant $\Delta$ of this system is
\[
\Delta = \Delta(\gamma_1,\ldots,\gamma_m) = \det(1 \ \gamma_i \ \gamma_i(\gamma_i+1) \ \ldots \ \gamma_i\cdots(\gamma_i+N-2))_{i=1,\ldots,N} = \prod_{1\leq i<j \leq N}(\gamma_j-\gamma_i) \neq 0.
\]
By choosing $b_N=a_0\neq 0$ we then obtain a unique solution, which we next construct explicitly as in \cite{V}.

For $\sigma=1,\ldots,N$, let $\Delta_\sigma(z)$ denote the determinant obtained from $\Delta$ after replacing $\gamma_\sigma$ by $z$. So
\[
\Delta_\sigma(z)=\Delta_{\sigma 0}+\Delta_{\sigma 1}z+\Delta_{\sigma 2}z(z+1)+\cdots+\Delta_{\sigma,N-1}z(z+1)\cdots(z+N-2),
\]
where $\Delta_{\sigma k}$ is the cofactor of $\Delta$ corresponding to the $\sigma,k$-entry \ $(\sigma=1,\ldots,N; k=0,1,\ldots,N-1)$. Clearly $\Delta_\sigma(\gamma_s) = 0$ for all $s \neq \sigma$, and therefore
\[
\Delta_\sigma(z) = c\prod_{s=1,s\neq\sigma}^N(z-\gamma_s),
\]
where
\[
c = \Delta \prod_{s=1,s\neq\sigma}^N(\gamma_\sigma-\gamma_s)^{-1},
\]
since $\Delta_\sigma(\gamma_\sigma) = \Delta$. This means that
\begin{equation}\label{7}
\Delta_{\sigma 0}+\Delta_{\sigma 1}z+\Delta_{\sigma 2}z(z+1)+\cdots+\Delta_{\sigma,N-1}z(z+1)\cdots(z+N-2) = \Delta \prod_{s=1,s\neq\sigma}^N\frac{z-\gamma_s}{\gamma_\sigma-\gamma_s}.
\end{equation}

We now choose $z = -\kappa$ in (\ref{7}) for each $\kappa = 0,1,\ldots,N-1$ to obtain
\[
\Delta_{\sigma 0}-\kappa\Delta_{\sigma 1}+\kappa(\kappa-1)\Delta_{\sigma 2}+\cdots+(-1)^\kappa \kappa!\Delta_{\sigma,\kappa} = \Delta \prod_{s=1,s\neq\sigma}^N\frac{\gamma_s+\kappa}{\gamma_s-\gamma_\sigma}.
\]
This gives
\begin{equation}\label{8}
A(\frac{\Delta_{\sigma0}}{\Delta},\frac{\Delta_{\sigma1}}{\Delta},\ldots,\frac{\Delta_{\sigma,N-1}}{\Delta})^T = 
\end{equation}
\[
(\frac{(-1)^0}{0!}\prod_{s=1,s\neq \sigma}^N \frac{\gamma_s}{\gamma_s - \gamma_\sigma},\frac{(-1)^1}{1!}\prod_{s=1,s\neq \sigma}^N \frac{ \gamma_s+1}{\gamma_s - \gamma_\sigma},\ldots,\frac{(-1)^{N-1}}{(N-1)!}\prod_{s=1,s\neq \sigma}^N \frac{\gamma_s+N-1}{\gamma_s - \gamma_\sigma})^T,
\]
where $A$ is the $N\times N$-matrix with rows
\[
(\frac{(-1)^\kappa}{\kappa!}, \frac{(-1)^{\kappa-1}}{(\kappa-1)!}, \cdots, \frac{(-1)^1}{1!}, \frac{(-1)^0}{0!}, 0,\ldots, 0), \ \kappa=0,1,\ldots,N-1.
\]
Since $A^{-1}$ is the matrix with rows
\[
(\frac{1}{k!},\frac{1}{(k-1)!}, \cdots, \frac{1}{1!}, \frac{1}{0!}, 0,\ldots, 0), \ k=0,1,\ldots,N-1,
\]
we immediately get
\begin{equation}\label{9}
\frac{k!\Delta_{\sigma k}}{\Delta} = \sum_{\tau = 0}^k(-1)^\tau{k \choose \tau}\prod_{s=1,s\neq \sigma}^N \frac{\gamma_s+\tau}{\gamma_s - \gamma_\sigma}, \quad k = 0,1,\ldots,N-1.
\end{equation}

By Cramer's rule (\ref{6}) has a solution
\[
b_k = -b_N\sum_{\sigma =1}^N\frac{\Delta_{\sigma k}}{\Delta}\prod_{\mu =0}^{N-1}(\gamma_\sigma + \mu), \ k = 0,1,\ldots,N-1.
\]
The choice $b_N=a_0=-1/N!$ and (\ref{9}) then give
\begin{equation}\label{10}
k!b_k = k!a_{N-k} = \sum_{\sigma =1}^N \sum_{\tau = 0}^k(-1)^\tau{k \choose \tau}\prod_{\mu =0}^{N-1}\frac{\gamma_\sigma + \mu}{1+\mu}\prod_{s=1,s\neq \sigma}^N \frac{\gamma_s+\tau}{\gamma_s - \gamma_\sigma}, \ k = 0,1,\ldots,N-1.
\end{equation}
Thus we have explicitly constructed Pad\'e approximations of the second kind $Q(z), P_1(z),\ldots, P_m(z)$ and the remainders
\[
R_j(z) = \sum_{\mu=N+n+1}^\infty c_{j\mu}z^\mu, \ j=1,\ldots,m.
\]
Note that here deg $Q(z) = N$, since if $a_N=b_0=0$, then $b_1,\ldots,b_N$ satisfy by (\ref{6}) a system of linear homogeneous equations with coefficient determinant $\Delta(\gamma_1+1,\ldots,\gamma_N+1) \neq 0$, and this implies $b_1=\cdots=b_N=0$, which contradicts our choice $b_N = 1/N!$.

This construction is not enough, since we shall need $m+1$ linearly independent approximations. To find further $m$ we fix $i, 1 \leq i \leq m$, and construct approximations $Q_i(z), P_{i1}(z),\ldots, P_{im}(z)$, where the gap in $Q_i(z)\varphi_i(z)$ is moved one step from the above construction and in $Q_i(z)\varphi_j(z), i \neq j$, it is in the same place as above. Thus we choose the coefficients $a_{ik}$ in
\[
Q_i(z) = \sum_{k=0}^N a_{ik}z^k
\]
in such a way that in
\[
Q_i(z)\varphi_j(z) = \sum_{\mu=0}^\infty c_{ij\mu}z^\mu
\]
$c_{ij\mu} = 0 \ (\mu=N+n-n_j+1+\delta_{ij},\ldots,N+n+\delta_{ij}; j=1,\ldots,m)$, here $\delta_{ij}$ is Kronecker's delta. Moreover, we assume an extra condition $c_{ii,N+1}=\alpha_i \neq 0$. By using the notations $a_{i,N-k} = b_k \ (k=0,1,\ldots,N)$ and
\[
\gamma_0 = n-n_i+1+\alpha_i,
\]
\[
\gamma_1 = n-n_1+1+\alpha_1+\delta_{1i}, \ldots, \gamma_{n_1} = n+\alpha_1+\delta_{1i},
\]
\[
\gamma_{n_1+1} = n-n_2+1+\alpha_2+\delta_{2i}, \ldots, \gamma_{n_1+n_2} = n+\alpha_2+\delta_{2i}, \ldots
\]
\[
\gamma_{n_1+\cdots+n_{m-1}+1} = n-n_m+1+\alpha_m+\delta_{mi}, \ldots, \gamma_N = n+\alpha_m+\delta_{mi},
\]
we see that $b_k$ should satisfy the system of equations
\[
b_0+b_1\gamma_0+b_2\gamma_0(\gamma_0+1)+\cdots+b_N\gamma_0(\gamma_0+1)\cdots(\gamma_0+N-1) = 1,
\]
\[
b_0+b_1\gamma_\sigma+b_2\gamma_\sigma(\gamma_\sigma+1)+\cdots+b_N\gamma_\sigma(\gamma_\sigma+1)\cdots(\gamma_\sigma+N-1) = 0, \ \sigma = 1,\ldots,N.
\]
The coefficient determinant $\Delta_i$ of this system is
\[
\Delta_i = \prod_{0\leq \ell<j\leq N} (\gamma_j-\gamma_{\ell}) \neq 0.
\]
Thus there exists a unique solution $b_0,b_1,\ldots,b_N$, which is obtained analogously to the above consideration leading to (\ref{9}).

Let $\Delta_i(z)$ be the determinant obtained from $\Delta_i$ after replacing $\gamma_0$ by $z$. Then
\[
\Delta_i(z)=\Delta_{i00}+\Delta_{i01}z+\Delta_{i02}z(z+1)+\cdots+\Delta_{i0N}z(z+1)\cdots(z+N-1),
\]
where $\Delta_{i0k}$ is the cofactor corresponding to the $0,k$-entry \ $(k=0,1,\ldots,N)$. Analogously to (\ref{7}) we now have
\[ 
\Delta_{i00}+\Delta_{i01}z+\Delta_{i02}z(z+1)+\cdots+\Delta_{i0N}z(z+1)\cdots(z+N-1) = \Delta_i \prod_{s=1}^N\frac{z-\gamma_s}{\gamma_0-\gamma_s}.
\]
As above, this gives
\[
\frac{k!\Delta_{i0k}}{\Delta_i} = \sum_{\tau = 0}^k(-1)^\tau{k \choose \tau}\prod_{s=1}^N \frac{\gamma_s+\tau}{\gamma_s - \gamma_0}, \quad k = 0,1,\ldots,N,
\]
and then, by Cramer's rule,
\begin{equation}\label{11}
k!b_k = k!a_{i,N-k} = \sum_{\tau = 0}^k(-1)^\tau{k \choose \tau}\prod_{s=1}^N \frac{\gamma_s+\tau}{\gamma_s - \gamma_0}, \quad k = 0,1,\ldots,N.
\end{equation}

So we have constructed $m$ systems of polynomials
\[
Q_i(z) = \sum_{k=0}^N a_{ik}z^k \neq 0, P_{ij}(z) = \sum_{\mu=0}^{N+n-n_j+\delta_{ij}} c_{ij\mu}z^\mu, \ j=1,\ldots,m; i=1,\ldots,m,
\]
such that deg $P_{ii}(z) = N+n-n_i+1, \deg P_{ij}(z) \leq N+n-n_j$, if $i \neq j$, and
\[
R_{ij}(z) := Q_i(z)\varphi_j(z)-P_{ij}(z) = \sum_{\mu=N+n+1+\delta_{ij}}^\infty c_{ij\mu}z^\mu.
\]
These systems with the first construction
\[
Q(z) =: Q_0(z), P_j(z) =: P_{0j}(z), R_j(z) =: R_{0j}(z), \ j=1,\ldots,m,
\]
satisfy the following independence lemma.\\

\noindent
\textbf{Lemma 1}. \textit{The determinant
\[
\Omega(z) = \det (Q_i(z) \ P_{i1}(z) \ \ldots \ P_{im}(z))_{i=0,1,\ldots,m} = a_N\alpha_1\cdots\alpha_mz^{m(N+n)+m},
\]
where $a_N \neq 0$ is given in (\ref{10}).}\\

\noindent
\textit{Proof}. Since $\deg Q_0(z) = N, \deg P_{0j} \leq N+n-n_j, \deg Q_i(z) \leq N, \deg P_{ii}(z) = N+n+n_i+1, \deg P_{ij} \leq N+n-n_j$, if $i \neq j \ (i,j = 1,\ldots,m)$, and the leading coefficients of $Q_0(z)$ and $P_{ii}(z)$ are $a_N$ and $c_{ii,N+n-n_i+1} = \alpha_i$, respectively, it follows that the leading term of $\Omega(z)$ is $a_N\alpha_1\cdots\alpha_mz^{m(N+n)+m} \neq 0$. On the other hand, $\Omega(z)$ can be given in $\mathbb{Q}[[z]]$ in the form
\[
\Omega(z) = (-1)^m\det (Q_i(z) \ R_{i1}(z) \ \ldots \ R_{im}(z))_{i=0,1,\ldots,m}.
\]
Since ord $R_{ij}(z) \geq N+n+1$, we have ord $\Omega(z) \geq m(N+n)+m$, which proves Lemma 1.

\section{Denominators}

In the proof of Theorems 1 and 2 we use the approximations of section 2, where $n_1=\cdots=n_m=n$, and so we consider in the following only this special case. For the study of the denominators of the coefficients of $Q_i(z)$ and $P_{ij}(z)$ we first recall a lemma from \cite [pp. 145-147]{M} considering the quotients
\[
\frac{(\alpha+1)_n}{n!} =: \frac{u_n}{v_n}, \ (u_n,v_n) = 1, v_n \geq 1, n = 0,1,\ldots,
\]
where $\alpha = r/s \neq -1,-2,\ldots$ with integers $r$ and $s \geq 1, (r,s) = 1$.\\

\noindent
\textbf{Lemma 2}. \textit{Let
\[
U_n = \prod_{p\nmid s}p^{[\log(\left|r\right|+sn)/\log p]}, \ V_n = s^{2n}.
\]
Then the least common multiples of $u_0,u_1,\ldots,u_n$ and of $v_0,v_1,\ldots,v_n$ are divisors of $U_n$ and $V_n$, respectively.}\\

For the following considerations we denote
\[
\alpha_j = \frac{r_j}{s_j}, (r_j,s_j) = 1, s_j \geq 1, \quad \alpha_k-\alpha_j = \frac{r_{kj}}{s_{kj}}, (r_{kj},s_{kj}) = 1, s_{kj} \geq 2.
\]
Further, let
\[
R = \max \{\left|r_j\right|\}, S = \max \{s_j\}, \quad \hat{R} = \max \{\left|r_{kj}\right|\}, \hat{S} = \max \{s_{kj}\}.
\]
Clearly $\hat{R} \leq 2RS$ and $\hat{S} \leq S^2$. 
  
Let us start by considering the coefficients of $Q_0(z)$ given in (\ref{10}). If $\gamma_\sigma = \kappa + \alpha_t, 1 \leq \kappa \leq n$, then
\begin{equation}\label{12}
\Pi_{\sigma\tau} := \prod_{s=1,s\neq \sigma}^N \frac{\gamma_s+\tau}{\gamma_s - \gamma_\sigma} = \prod_{j=1,j\neq t}^m\prod_{\nu=1}^n\frac{\alpha_j+\tau+\nu}{\alpha_j-\alpha_t+\nu-\kappa}\cdot \prod_{\nu=1,\nu\neq\kappa}^n\frac{\alpha_t+\tau+\nu}{\nu-\kappa} =
\end{equation}
\[
\prod_{j=1,j\neq t}^m\frac{(\alpha_j+\tau+1)_n}{(\alpha_j-\alpha_t-\kappa+1)_n}\cdot (-1)^{\kappa-1}\frac{\kappa}{\alpha_t+\tau+\kappa}{n \choose \kappa}\frac{(\alpha_t+\tau+1)_n}{n!}.
\]
Since all the numbers $r_t+s_t(\tau+\kappa)$ are factors of
\[
\prod_p p^{[\log(R+S(N+n))/\log p]},
\]
it follows by Lemma 2 and (\ref{10}), that
\[
k!D_1a_{N-k} \in \mathbb{Z}, \ k = 0,1,\ldots, N,
\]
where
\begin{equation}\label{13}
D_1 = s^{2N}(s_1\cdots s_m)^{2n}\prod_p p^{(m-1)[\log(\hat{R}+\hat{S}n)/\log p]+[\log(R+S(N+n))/\log p]},
\end{equation}
remember that $s = lcm \{s_1,\ldots,s_m\}$. Thus
\[
N!D_1Q_0(z), \ N!D_1s^NP_{0j}(z) \in \mathbb{Z}, \ j = 1,\ldots, m.
\]

For $1 \leq i \leq m$ the coefficients of $Q_i(z)$ are given in (\ref{11}), where
\begin{equation}\label{14}
\prod_{s=1}^N \frac{\gamma_s+\tau}{\gamma_s - \gamma_0} = \prod_{j=1,j\neq i}^m\prod_{\nu=1}^n\frac{\alpha_j+\tau+\nu}{\alpha_j-\alpha_i+\nu-1}\cdot\prod_{\nu=2}^{n+1}\frac{\alpha_i+\tau+\nu}{\nu-1} = \prod_{j=1,j\neq i}^m\frac{(\alpha_j+\tau+1)_n}{(\alpha_j-\alpha_i)_n}\cdot\frac{\alpha_i+\tau+2)_n}{n!}.
\end{equation}
The above consideration together with (\ref{11}) and our choice $c_{ii,N+1} = \alpha_i \ (i=1,\ldots,m)$ then imply the following lemma.\\

\noindent
\textbf{Lemma 3}. \textit{We have
\[
k!D_1a_{i,N-k} \in \mathbb{Z}, \ k=0,1,\ldots,N; i=0,1,\ldots,m.
\]
Further, if $D_2 := D_1s^N$, then}
\[
N!D_2Q_i(z), \ N!D_2P_{ij}(z) \in \mathbb{Z}[z], \ i = 0,1,\ldots,m; j = 1,\ldots, m. \\
\]

\section{Upper bounds}

We shall first obtain upper bounds for $\left|Q_i(z)\right|$ and $\left|P_{ij}(z)\right|$. By Lemma 2 and (\ref{12}),
\[
\left|\Pi_{\sigma\tau}\right| \leq \prod_{j=1,j\neq t}^m(\frac{s_{jt}^{2n}}{s_j^{n}}\prod_{p} p^{[\log(R+S(N+n))/\log p]})\cdot \frac{n2^n}{s_t^{n}}\prod_{p} p^{[\log(R+S(N+n))/\log p]} \leq 
\]
\[
(s_1\cdots s_m)^ns_t^{2(m-2)n}n2^n \prod_{p} p^{m[\log(R+S(N+n))/\log p]}.
\]
Therefore (\ref{10}) and one further application of Lemma 2 imply
\[
(N-k)!\left|a_k\right| \leq nN^22^{N+n}S^{(2m-4)n}\prod_{p} p^{(m+1)[\log(R+S(N+n))/\log p]}, \ k=0,1,\ldots,N.
\]
Further,
\[
\left|c_{j\mu}\right| = \left|\sum_{k=0}^\mu (\mu-k)!a_k\frac{(\alpha_j)_{\mu-k}}{(\mu-k)!}\right| \leq n(1+N)N^22^{N+n}S^{(2m-4)n}\cdot
\]
\[
\prod_{p} p^{(m+1)[\log(R+S(N+n))/\log p]+[\log (R+SN)/\log p]}, \ \mu = 0,1,\ldots,N; j = 1,\ldots,m. 
\]

In the case $1 \leq i \leq m$ a similar consideration using (\ref{11}) and (\ref{14}) gives
\[
\left|(N-k)!a_{ik}\right| \leq (1+N)2^NS^{(2m-4)n}\prod_{p} p^{(m-1)[\log(R+S(N+n))/\log p]+[\log (R+S(N+n+1)/\log p]}, \ k=0,1,\ldots,N,
\]
and
\[
\left|c_{ij\mu}\right| \leq (1+N)^22^NS^{(2m-4)n}\cdot
\]
\[
\prod_{p} p^{(m-1)[\log(R+S(N+n))/\log p]+[\log (R+S(N+n+1)/\log p]+[\log (R+SN)/\log p]}, \ \mu = 0,1,\ldots,N+\delta_{ij}; j = 1,\ldots,m.
\]

By the weak form of the prime number theorem, see for example \cite[p. 296]{Bu}, the number of primes $p \leq x$
\[
\pi(x) \leq 8\log 2\frac{x}{\log x} < 6\frac{x}{\log x}
\]
for all $x>1$, and therefore the above estimates imply immediately the following lemma.\\

\noindent
\textbf{Lemma 4}. \textit{We have, for all $i=0,1,\ldots,m; j=1,\ldots,m$,
\[
\left|Q_i(z)\right|, \left|P_{ij}(z)\right| \leq c_1n^5c_2^n\max \{1,\left|z\right|^{mn+1}\},
\]
where}
\[
c_1 = m(1+m)^2(2+m)e^{6(m+2)R+6S}, c_2 = 2^{m+1}S^{2m-4}e^{6(m^2+3m+1)S}.\\
\]

To consider $\left|R_{ij}(z)\right|_p$ we note that, for all $\mu \geq N+n+1$,
\[
c_{ij\mu} = \sum_{k=0}^N a_{ik}(\alpha_j)_{\mu-k} = (\mu-N)!\sum_{k=0}^N (N-k)!a_{ik}\frac{(\mu-N+1)\cdots(\mu-k)}{(N-k)!}\frac{(\alpha_j)_{\mu-k}}{(\mu-k)!}.
\]
Therefore, by Lemma 3, $\left|D_1c_{ij\mu}\right|_p \leq \left|(\mu-N)!\right|_p$, if $p \nmid s$. Thus, for $p \nmid s, \left|z\right|_p \leq 1$,
\[
\left|D_1R_{ij}(z)\right|_p \leq \max \{\left|(\mu-N)!z^\mu\right|_p \mid \mu \geq N+n+1\} = \left|(n+1)!z^{N+n+1}\right|_p
\]
giving the following result.\\

\noindent
\textbf{Lemma 5}. \textit{If $p \nmid s$ and $\left|z\right|_p \leq 1$, then}
\[
\left|D_1R_{ij}(z)\right|_p \leq \left|(n+1)!\right|_p\left|z\right|_p^{mn+n+1}.\\
\]

\section{Proof of Theorems 1 and 2}

Let $a\in\mathbb{Z}\setminus\{0\}$ and denote
\[
Q_i = N!D_2Q_i(a), P_{ij} = N!D_2P_{ij}(a), R_{ij} = N!D_2R_{ij}(a), \quad i=0,1,\ldots,m; j=1,\ldots,m.
\]  
By Lemmas 1 and 3 the numbers $Q_i$ and $P_{ij}$ are integers satisfying
\begin{equation}\label{15}
\det (Q_i \ P_{i1} \ \ldots \ P_{im})_{i=0,1,\ldots,m} \neq 0,
\end{equation}
and $R_{ij}$ are defined for all $p \nmid s$, by Lemma 5.

We now assume that the assumptions of Theorem 1 are valid and a linear form
\[
L = \ell_0 + \ell_1\varphi_1(a) + \cdots + \ell_m\varphi_m(a)
\]
with $(\ell_0,\ell_1,\ldots,\ell_m)\in\mathbb{Z}^{m+1}\setminus\{\underline{0}\}$ satisfies $L_p = 0$ for all $p\in P_1$. Since
\[
Q_iL_p = Q_i\ell_0 + \sum_{j=1}^mP_{ij}\ell_j + \sum_{j=1}^mR_{ij}\ell_j =: \Lambda + \sum_{j=1}^mR_{ij}\ell_j,
\]
the above assumption means that $\Lambda$ is an integer and
\[
\Lambda = -\sum_{j=1}^mR_{ij}\ell_j, \quad p\in P_1.
\]

By (\ref{15}) there exists an $i, 0 \leq i \leq m$, such that $\Lambda \neq 0$. Thus, by Lemmas 4 and 5,
\[
1 = \left|\Lambda\right|\prod_p\left|\Lambda\right|_p \leq \left|\Lambda\right|\prod_{p\in P_1}\left|\Lambda\right|_p \leq (1+m)c_1hn^5c_2^n\left|a\right|^{mn+1}(mn)!D_2\prod_{p\in P_1}\left|a^{mn+n+1}(mn)!(n+1)!\right|_p,
\]
where $h = \max \{\left|\ell_0\right|,\left|\ell_1\right|,\ldots,\left|\ell_m\right|\}$. Here
\[
D_2 = D_1s^N \leq s^{3N}S^{2mn}e^{6((m-1)\hat{R}+R+(m-1)\hat{S}n+S(N+n))}
\]
by (\ref{13}). Therefore we have
\begin{equation}\label{16}
1 \leq c_3hc^n\left|a\right|^{mn+1}(mn)!\prod_{p\in P_1}\left|a^{mn+n+1}(mn)!(n+1)!\right|_p
\end{equation}
for all $n \geq 5$, where $c_3 = c_1(1+m)e^{12(m-1)RS+6R}$ and
\begin{equation}\label{17}
c = c_2s^{3m}S^{2m}e^{2+6S((m-1)S+m+1)} = 2^{m+1}s^{3m}S^{4m-4}e^{2+6S(m^2+2m+2+(m-1)S)}.
\end{equation}
The above inequality and (\ref{4}) give a contradiction for large $n$ proving Theorem 1.

To prove Theorem 2 we fix $n_1$ to be an integer satisfying
\begin{equation}\label{18}
c_4hc_5^{n_1}e^{\frac{5}{4}(m+1)n_1\sqrt{\log n_1}-(n_1+1)\log (n_1+1)} < 1, \quad c_4 = 6e^{10}c_3\left|a\right|, c_5 = 2ce^{9(m+1)+1}s^m\left|a\right|^m,
\end{equation}
and repeat the above consideration with $n=n_1$ and $P_1 = \{p\mid e^{\sqrt{\log n}} < p \leq mn, p \nmid s\}$. Then we get again (\ref{16}) with this $P_1$. Now
\[
\prod_{p\mid s}\left|(mn)!(n+1)!\right|_p^{-1} \leq \prod_{p\mid s}p^{(mn+n+1)/(p-1)} \leq 2^{n+1}3s^{mn}
\]
and
\[
\prod_{p\leq e^{\sqrt{\log n}}}\left|(mn)!(n+1)!\right|_p^{-1} \leq \prod_{p\leq e^{\sqrt{\log n}}}e^{(mn+n+1)\frac{\log p}{p-1}} \leq e^{(mn+n+1)\sum_{p\leq e^{\sqrt{\log n}}}\frac{\log p}{p-1}} \leq e^{(mn+n+1)(1+\frac{5}{4}\sum_{p\leq e^{\sqrt{\log n}}}\frac{\log p}{p})}.
\]
Since, for a positive integer $X$,
\[
\log X! = \sum_{p\leq X}\log p\sum_{j\geq 1}\left[\frac{X}{p^j}\right] \geq \sum_{p\leq X}\left[\frac{X}{p}\right]\log p \geq X\sum_{p\leq X}\frac{\log p}{p} - \pi(X)\log X,
\]
the use of the above mentioned weak form of the prime number theorem together with Stirling's formula implies
\[
X\sum_{p\leq X}\frac{\log p}{p} \leq X\log X - X + \log \sqrt{2\pi X} + \frac{1}{12X} + 6X \leq X\log X + 6X
\]
for all $X \geq 2$. Thus, if $x\geq 2$, then
\[
\sum_{p\leq x}\frac{\log p}{p} = \sum_{p\leq \left[x\right]}\frac{\log p}{p} \leq \log x + 6,
\]
and so it follows, by (\ref{16}), that
\[
1 \leq 6e^9c_3\left|a\right|h(2ce^{9(m+1)}s^m\left|a\right|^m)^ne^{\frac{5}{4}(m+1)n\sqrt{\log n}}\frac{1}{(n+1)!}.
\]
A further application of Stirling's formula then gives
\begin{equation}\label{19}
1 \leq c_4hc_5^ne^{\frac{5}{4}(m+1)n\sqrt{\log n}-(n+1)\log (n+1)},
\end{equation}
which is a contradiction with our choice of $n=n_1$ above. Therefore $L_p \neq 0$ for some $p, e^{\sqrt{\log n_1}} < p \leq mn_1$.

Let now $n_2$ be an integer such that $mn_1 < e^{\sqrt{\log n_2}}$ and (\ref{18}) holds if $n_1$ is replaced by $n_2$ . As above we see that at least one $p, e^{\sqrt{\log n_2}} < p < mn_2$, satisfies $L_p \neq 0$. Continuing in this way we obtain an infinite sequence of primes $p_1 < p_2 <\cdots$ with $L_{p_i} \neq 0$.

To prove the second claim of Theorem 2 we take  
\[
n = \left[\frac{\log H}{\log\log H}(1+\epsilon(H))\right], \ H = \max \{h,H_0\},
\]
where $H_0$ and the function $\epsilon(H) > 0$ satisfying $\lim_{H\rightarrow\infty}\epsilon(H) = 0$ are specified later in such a way that $n$ satisfies the conditions of $n_1$ above. Let now $P_1 = \{p\mid e^{\sqrt{\log n}} < p \leq mn, p \nmid s\}$. Similarly as above the inequality (\ref{19}) holds.

Assuming that $0 < \epsilon(H) \leq 1$ we now have
\[
\log c_4 + \log h + n\log c_5 + \frac{5}{4}(m+1)n\sqrt{\log n} - (n+1)\log (n+1) \leq \log c_4 + \log h + 2\log c_5 \left(\frac{\log H}{\log\log H}\right) +
\]
\[
\frac{5}{2}(m+1)\frac{\log H}{\log\log H}\sqrt{\log\log H} - \frac{\log H}{\log\log H}(1 + \epsilon(H))(\log\log H - \log\log\log H + \log (1 + \epsilon(H))) \leq 
\]
\[
\log c_4 + 2\log c_5 \frac{\log H}{\log\log H} + \frac{5}{2}(m+1)\frac{\log H}{\sqrt{\log\log H}}  + 2\frac{\log H \log\log\log H}{\log\log H} - \epsilon(H)\log H < 0,
\]
if we choose
\[
\epsilon(H) = \frac{\frac{5}{2}(m+1)+6}{\sqrt{\log\log H}}
\]
and
\begin{equation}\label{20}
\log H_0 = \max \{(\log c_4)^2, e^{(\log c_5)^2}, e^{(\frac{5}{2}(m+1)+6)^2}\}.
\end{equation}
This gives a contradiction with (\ref{19}), and so there must be a prime $p\in P_1$ with $L_p \neq 0$. This proves Theorem 2.

\begin{small}

\vskip0.6cm
\noindent
Keijo V\"a\"an\"anen \\
Department of Mathematical Sciences  \\
University of Oulu \\
P. O. Box 3000 \\
90014 Oulun yliopisto, Finland \\
E-mail: keijo.vaananen@oulu.fi \\
\end{small}


\begin{center}
\textbf{References}
\end{center}
\vspace{-1.5cm}
\begin{thebibliography}{}
\bibitem{BCY} D. Bertrand, V. Chirskii and J Yebbou, \textit{Effective estimates for global relations on Euler-type series}, Ann. Fac. Sci. Toulouse Math. (6) 13 (2004), 241-260. 

\bibitem{Bu} P. Bundschuh, \textit{Einf\"uhrung in die Zahlentheorie}, 4 Aufl., Springer-Lehrbuch, Springer, 1998.

\bibitem{C1} V. G. Chirskii, \textit{On the arithmetic properties of generalized hypergeometric series with irrational parameters}, Izvestiya: Mathematics 78:6 (2014), 1244-1260.

\bibitem{C2} V. G. Chirskii, \textit{Arithmetic properties of Euler series}, Moscov. Univ. Math. Bull. 70, no. 1 (2015), 41-43.

\bibitem{F} N. I. Fel'dman, \textit{Lower estimates for some linear forms}, Vestnik Moscov. Univ. Ser. I, Mat. Mekh. 22, no. 2 (1967), 63-72.

\bibitem{M} K. Mahler, \textit{Lectures on Transcendental Numbers}, Lecture Notes in Mathematics 546, Springer, 1976.

\bibitem{MZ} T. Matala-aho and W. Zudilin, \textit{Euler's factorial series and global relations}, arXiv:1703.02633v1 (2017).

\bibitem{V} K. V\"a\"an\"anen, \textit{On a result of Fel'dman on linear forms in the values of some $E$-functions}, arXiv:1704.01762v1 (2017).


\end{thebibliography}
\end{document}